\documentclass{amsart}
\usepackage{mathrsfs}
\usepackage{lamsarrow}
\usepackage{pb-diagram}
\usepackage{amssymb}
\usepackage{pb-lams}

\newtheorem{theorem}{Theorem}[section]
\newtheorem{lemma}[theorem]{Lemma}

\newtheorem{definition}[theorem]{Definition}
\newtheorem{remark}[theorem]{Remark}

\newtheorem{construction}[theorem]{Construction}

\DeclareMathOperator{\SW}{\mathcal {SW}}
\DeclareMathOperator{\Hom}{Hom}
\DeclareMathOperator{\RHom}{RHom}

\DeclareMathOperator{\map}{map}
\DeclareMathOperator{\Ho}{Ho}
\DeclareMathOperator{\HoC}{Ho(\C)}

\DeclareMathOperator{\colim}{colim}

\newcommand{\mZ}{{\mathbb Z}}
\newcommand{\mS}{{\mathbb S}}
\newcommand{\C}{{\mathcal C}}
\newcommand{\D}{{\mathcal D}}
\newcommand{\T}{{\mathcal T}}

\newcommand{\Spc}{{\mathcal S}p}
\newcommand{\sset}{\mathcal{S}_*}
\newcommand{\iso}{\cong}
\newcommand{\rscript}[1]{^{\mbox{\scriptsize {#1}}}}

\newcommand{\sm}{\wedge}

\newcommand{\tensor}{\otimes}
\newcommand{\boxprod}{\mathbin{\square }}
\newcommand{\Wedge}{{\scriptstyle \vee}}

\newcommand{\mc}{\colon \,}

\renewcommand{\to}{\longrightarrow}
\newcommand{\varrow}[1]{\hbox to #1{\rightarrowfill}}
\newcommand{\varl}[2]{\stackrel{#2}{\hbox to #1{\leftarrowfill}}}
\newcommand{\varr}[2]{\stackrel{#2}{\hbox to #1{\rightarrowfill}}}
\newcommand{\parallelarrows}[1]{\begin{array}{c} {\hbox to
#1{\rightarrowfill}}  \vspace{-0.35cm} \\ {\hbox to
#1{\rightarrowfill}} \end{array}}

\begin{document}

\title{A uniqueness theorem for stable homotopy theory}
\date{\today; 1991 AMS Math.\ Subj.\ Class.: 55U35, 55P42}
\author{Stefan Schwede}
\thanks{Research supported by a BASF-Forschungsstipendium der
Studienstiftung des deutschen Volkes}
\address{Fakult\"at f\"ur Mathematik \\ Universit\"at Bielefeld\\
33615 Bielefeld, Germany}
\email{schwede@mathematik.uni-bielefeld.de}
\author{Brooke Shipley}
\thanks{Research partially supported by NSF grants}
\address{Department of Mathematics \\ University of Chicago\\  
Chicago, IL 60637 
\\ USA}
\email{bshipley@math.uchicago.edu}
\maketitle

\vspace*{-0.5cm}

\section{Introduction}

Roughly speaking, the stable homotopy category is
obtained from the homotopy category of topological spaces by inverting 
the suspension functor, yielding a `linear' approximation to 
the homotopy category of spaces.  The isomorphism classes of objects in
the stable homotopy category represent the generalized cohomology theories, 
defined by the Eilenberg-Steenrod axioms~\cite{e-s} without the dimension
axiom (which distinguishes `ordinary' from `generalized' cohomology theories).

The first construction of the full stable homotopy category was given 
by Boardman~\cite{boardman-spectra}. 
More recently many models for the stable homotopy category have been found, 
most of which have the additional structure of a {\em closed model category}
in the sense of Quillen~\cite{Q}.
In~\cite{margolis}, H.~R.~Margolis introduced a short list of axioms
and conjectured that they characterize the stable homotopy category
up to an equivalence of categories.    
In Theorem~\ref{thm-Margolis} we prove that these axioms indeed
uniquely specify the stable homotopy category whenever
there is some underlying Quillen model category.

We also prove a more structured version, the Uniqueness Theorem below, 
which states that the model category of spectra 
itself is uniquely determined by certain equivalent conditions, 
up to so called {\em Quillen equivalence} of model categories
(a particular adjoint pair of functors which induces 
equivalences of homotopy categories, see Definition~\ref{def-Quillen pair}).
This implies that the higher order structure
including cofibration and fibration sequences, homotopy (co-)limits,
Toda brackets, and mapping spaces is uniquely determined
by these conditions.
This is of interest due to the recent plethora of new model categories 
of spectra~\cite{ekmm, hss, mmss, lydakis-simplicial}. 
The Uniqueness Theorem provides criteria on the homotopy category level
for deciding whether a model category captures the stable homotopy theory
of spectra; the search for such intrinsic characterizations was another 
main motivation for this project.  

A model category is {\em stable} if the suspension functor is invertible
up to homotopy. For stable model categories the homotopy category
is naturally triangulated and comes with an action by the graded
ring $\pi^s_{\ast}$ of stable homotopy groups of spheres,
see \ref{pi-structure on Ho(C)}.  The Uniqueness Theorem shows that this
$\pi^s_{\ast}$-triangulation determines the stable homotopy theory up 
to Quillen equivalence.   

\medskip

\noindent{\bf Uniqueness Theorem.} {\em Let $\C$ be a 
stable model category. Then the following four conditions are equivalent:
\renewcommand{\labelenumi}{(\arabic{enumi})}
\begin{enumerate}
\item There is a chain of Quillen equivalences between $\C$ 
and the model category of spectra.
\item There exists a $\pi_{\ast}^s$-linear equivalence 
between the homotopy category of $\C$ and the homotopy category of spectra.
\item The homotopy category of $\C$ has a small weak generator $X$ for which 
$[X,X]_{\ast}^{\HoC}$ is freely generated as a 
$\pi_{\ast}^s$-module by the identity map of $X$.
\item  The model category $\C$ has a cofibrant-fibrant small weak generator $X$
for which the unit map $\mS\to \Hom(X,X)$ is a $\pi_*$-isomorphism 
of spectra. 
\end{enumerate}
}

\bigskip

The Uniqueness Theorem is proved in a slightly more general form as 
Theorem~\ref{thm-recognition}. The extra generality consists of a local 
version at subrings of the ring of rational numbers.  
Our reference model for the category of spectra 
is that of Bousfield and Friedlander \cite[Def.\ 2.1]{BF};
this is probably the simplest model category of spectra 
and we review it in Section \ref{sec-local}.
There we also discuss the {\em $R$-local} model structure 
for spectra for a subring $R$ of the ring of rational numbers, 
see Lemma \ref{lem-R-local model structure}.
The notions of `smallness' and `weak generator' are recalled 
in \ref{def-generator, small, margolis}.
The unit map is defined in \ref{def-unit}.

Our work here grows out of recent developments in axiomatic stable homotopy 
theory. Margolis' axiomatic approach was generalized in~\cite{hps} 
to study categories which share the main formal properties 
of the stable homotopy category, namely triangulated symmetric monoidal 
categories with a weak generator or a set of weak generators.  
Hovey~\cite[Ch.\ 7]{hovey-book} then studied properties of model categories 
whose homotopy categories satisfied these axioms.  
Heller has given an axiomatization of the concept of a 
``homotopy theory'' \cite{heller-memoirs}, and then characterized 
the passage to spectra by a universal property in his context, 
see~\cite[Sec.\ 8-10]{heller}. The reader may want to compare this with
the universal property of the {\em model category} of spectra which we prove
in Theorem \ref{thm-universal} below.

In~\cite{ss-class} we classify stable model categories with 
a small weak generator as modules over a ring spectrum, 
see Remark~\ref{rem-ss-class}.  
The equivalence of parts (1) and (4) of our 
Uniqueness Theorem here can be seen as a special case of this classification.
In general, though, there are no 
analogues of parts (2) and (3) of the Uniqueness Theorem. 
Note, that here, as in~\cite{ss-class}, we ignore the smash product
in the stable homotopy category; 
several comparisons and classification results respecting smash products 
can be found in~\cite{sch-comparison, mmss}.

\section{Stable model categories} \label{sec-stable}

Recall from \cite[I.2]{Q} or \cite[6.1.1]{hovey-book}
that the homotopy category of a pointed
model category supports a suspension functor $\Sigma$ with a
right adjoint loop functor $\Omega$.

\begin{definition} \label{def-stable model category}
{\em  A {\em stable model category} is a pointed, complete and cocomplete 
category with a model category structure for which
the functors $\Omega$ and $\Sigma$ on the homotopy category are
inverse equivalences.}
\end{definition}

The homotopy category of a stable model category has  
a large amount of extra structure,
some of which will play a role in this paper. 
First of all, it is naturally a triangulated category (cf.\ \cite{verdier}). 
A complete reference for this fact can be found in \cite[7.1.6]{hovey-book}; 
we sketch the constructions: by definition the suspension 
functor is a self-equivalence of the
homotopy category and it defines the shift functor. 
Since every object is a two-fold suspension,
hence an abelian co-group object, the homotopy category of 
a stable model category is additive.
Furthermore, by \cite[7.1.11]{hovey-book} the cofiber sequences and 
fiber sequences of \cite[I.3]{Q} coincide
up to sign in the stable case, and they define the distinguished triangles.
Since we required a stable model category
to have all limits and colimits, its homotopy category will have infinite sums 
and products.

Apart from being triangulated, the homotopy category of a stable model
category has a natural action of the ring $\pi_*^s$ 
of stable homotopy groups of spheres. Since this action is central to
this paper, we formalize and discuss it in some detail.
For definiteness we set $\pi_n^s=\colim_k \, [S^{n+k},S^k]$,
where the colimit is formed along right suspension
$-\sm 1_{S^1}:[S^{n+k},S^k]\to [S^{n+k+1},S^{k+1}]$. 
The ring structure is given by composition of representatives.

\begin{definition} \label{def-pi_*}
{\em A {\em $\pi_*^s$-triangulated category} is
a triangulated category $\T$ with bilinear pairings
\[ \pi_n^s \ \tensor \ \T(X,Y)  \ \varrow{1cm} \ \T(X[n],Y) \ , 
\quad \alpha\tensor f \longmapsto \alpha\cdot f \]
for all $X$ and $Y$ in $\T$ and all $n\geq 0$, where $X[n]$ is the $n$-fold
shift of $X$. Furthermore the pairing must be associative, unital
and composition in $\T$ must be bilinear in the sense that 
$(\alpha\cdot g) \circ f[n] = \alpha\cdot(g\circ f) = g\circ(\alpha\cdot f)$.
A {\em $\pi_*^s$-exact functor} between  $\pi_*^s$-triangulated categories
is a functor $L:\T\to \overline{\T}$ together with a natural isomorphism
$\tau: L(X)[1]\iso L(X[1])$ such that
\begin{itemize}
\item $(L,\tau)$ forms an exact functor of triangulated categories, i.e.,
for every distinguished triangle $X\to Y\to Z\to X[1]$ in $\T$
the sequence $L(X)\to L(Y) \to L(Z) \to L(X)[1]$ is a distinguished
triangle in $\overline{\T}$, where the third map is the composite 
$L(Z)\to L(X[1]) \stackrel{\tau^{-1}}{\to} L(X)[1]$;
\item $(L,\tau)$ is $\pi_*^s$-linear, i.e., for all $X$ and $Y$ in $\T$ and $n\geq 0$ the following diagram commutes
\[\begin{diagram}
\node{\pi_n^s\tensor\T(X,Y)\hspace*{1.05cm}}  \arrow[2]{e,t}{\cdot}
\arrow{s,l}{\mbox{\scriptsize Id}\,\tensor\, L} 
\node[2]{\T(X[n],Y)}
\arrow{s,r}{L} \\
\node{\pi_n^s\tensor\overline{\T}(L(X),L(Y))} 
\arrow{e,b}{\cdot} 
\node{\overline{\T}(L(X)[n],L(Y))} 
\node{\overline{\T}(L(X[n]),L(Y))}  \arrow{w,b}{\circ\tau}
\end{diagram}\]
where $\tau:L(X)[n]\to L(X[n])$ is the $n$-fold iterate
of instances of the isomorphism $\tau$.
\end{itemize} 
A {\em $\pi_*^s$-linear equivalence} between  $\pi_*^s$-triangulated categories
is a $\pi_*^s$-exact functor which is an equivalence of categories and
whose inverse is also $\pi_*^s$-exact.}
\end{definition}

Every triangulated category can be made into a $\pi_*^s$-triangulated category
in a trivial way by setting $\alpha\cdot f=0$ whenever the dimension 
of $\alpha$ is positive.
Now we explain how the homotopy category of a stable model category
is naturally a $\pi_*^s$-triangulated category.

\begin{construction} \label{pi-structure on Ho(C)}
{\em Using the technique of {\em framings}, Hovey \cite[5.7.3]{hovey-book}
constructs a pairing
\[ \sm^L \ : \ \HoC \ \times \ \Ho(\sset) \ \varrow{1cm} \ \HoC \]
which makes the homotopy category of a pointed model category $\C$ 
into a module (in the sense of \cite[4.1.6]{hovey-book}) 
over the symmetric monoidal
homotopy category of pointed simplicial sets under smash product.
In particular, the pairing is associative and unital up to coherent
natural isomorphism, and smashing with the simplicial circle $S^1$ is
naturally isomorphic to suspension as defined by Quillen \cite[I.2]{Q}.
If $\C$ is stable, we may take $X[1]:= X\sm^L S^1$ as the shift functor 
of the triangulated structure.
We define the action 
\[ \pi_n^s \ \tensor \ [X,Y]^{\HoC} \ \  \quad \varr{1.5cm}{} \quad
[X[n],Y]^{\HoC} \]
as follows.
Suppose $\alpha:S^{n+k}\to S^k$ is a morphism in the homotopy category 
of pointed simplicial sets which represents an element of  
$\pi_n^s=\colim_k \, [S^{n+k},S^k]$ 
and $f:X\to Y$ is a morphism in the homotopy category of $\C$.
Since $\C$ is stable, smashing with $S^k$ is a bijection of 
morphism groups in the homotopy category. 
So we can define $\alpha\cdot f$ to be the unique morphism 
in $[X \sm^L S^n,Y]^{\HoC}$ such that
$(\alpha\cdot f)\sm^L 1_{S^k}= f\sm^L\alpha$ in 
the group $[X \sm^L S^{n+k},Y\sm^L S^k]^{\HoC}$.
Here and in the following we identify the $n$-fold shift
$X[n]=(\cdots((X\sm^L S^1)\sm^L S^1)\cdots )\sm^L S^1$ with
$X\sm^L S^n$  under the associativity isomorphism which is constructed 
in the proof of \cite[5.5.3]{hovey-book} 
(or rather its pointed analog \cite[5.7.3]{hovey-book}); 
this way we regard $\alpha\cdot f$ as an element of the group 
$[X[n],Y]^{\HoC}$. Observe that $\pi_*^s$ acts from the {\em left}
even though simplicial sets act from the {\em right} on the
homotopy category of $\C$.

By construction $\alpha\cdot f=(\alpha\sm 1_{S^1})\cdot f$,
so the morphism  $\alpha\cdot f$ only depends on the class of $\alpha$ in the
stable homotopy group $\pi_n^s$.
The $\pi_*^s$-action is unital; associativity can be seen as follows:
if  $\beta\in[S^{m+n+k},S^{n+k}]$ represents another stable homotopy element,
then we have
\begin{eqnarray} (\beta\cdot f)\sm^L \alpha & = & 
(1_Y \sm^L \alpha) \circ ((\beta\cdot f)\sm^L 1_{S^{n+k}}) \quad = \quad
(1_Y \sm^L \alpha) \circ (f\sm^L \beta) \nonumber \\
 & = & f \sm^L (\alpha\circ\beta)
\quad = \quad ((\alpha\circ\beta)\cdot f) \sm^L 1_{S^k} \nonumber
\end{eqnarray}
in the group $[X\sm^L S^{m+n+k},Y\sm^L S^k]^{\HoC}$.
According to the definition of $\alpha\cdot(\beta\cdot f)$ this means
that  $\alpha\cdot(\beta\cdot f)=(\alpha\circ\beta)\cdot f$.
Bilinearity of the action is proved in a similar way.}
\end{construction}

\begin{definition} \label{def-Quillen pair}
{\em A pair of adjoint functors between model categories is a 
{\em Quillen adjoint pair} if the right adjoint preserves fibrations and
trivial fibrations.
An equivalent condition is that the left adjoint preserves cofibrations 
and trivial cofibrations.   
A Quillen adjoint pair induces an adjoint pair of functors 
between the homotopy categories~\cite[I.4 Thm.\ 3]{Q}, 
the {\em total derived functors}.
A Quillen functor pair is a {\em Quillen equivalence} if the total 
derived functors are adjoint equivalences of the homotopy categories.}
\end{definition}
 
The definition of Quillen equivalences just given is not the most common one;
however it is equivalent to the usual definition by~\cite[1.3.13]{hovey-book}.
Suppose $F:\C\to\D$ is the left adjoint of a Quillen adjoint pair between
pointed model categories.
Then the total left derived functor $LF:\HoC\to\Ho(\D)$ of $F$ 
comes with a natural isomorphism $\tau:LF(X)\sm^L S^1\to LF(X\sm^L S^1)$ 
with respect to which it preserves cofibration sequences, 
see \cite[I.4 Prop.\ 2]{Q} or \cite[6.4.1]{hovey-book}. 
If $\C$ and $\D$ are stable, this makes $LF$ into an exact functor 
with respect to $\tau$.
It should not be surprising that $(LF,\tau)$ is also $\pi_*^s$-linear 
in the sense of Definition \ref{def-pi_*}, but showing this requires a careful
review of the definitions which we carry out in Lemma \ref{left derived}.

\begin{remark}\label{morava-non-example}
{\em 
In Theorem~\ref{thm-recognition} below we 
show that the $\pi_{\ast}^s$-triangulated homotopy category determines
the Quillen equivalence type of the model category of spectra.  
This is not true for general stable model categories.  
As an example we consider the $n$-th Morava $K$-theory spectrum $K(n)$ 
for $n>0$ and some fixed prime $p$.
This spectrum admits the structure of an $A_{\infty}$-ring spectrum
\cite{robinson-K(n)}, and so its module spectra form a stable model category.
The coefficient ring $K(n)_{\ast}={\mathbb F}_p[v_n,v_n^{-1}]$, 
with $v_n$ of degree $2p^n-2$, is a graded field, and so the homotopy category
of $K(n)$-modules is equivalent, via the homotopy group functor, 
to the category of graded  $K(n)_{\ast}$-modules.
Similarly the derived category of differential graded  $K(n)_{\ast}$-modules
is equivalent, via the homology functor, to the category of graded  
$K(n)_{\ast}$-modules.
This derived category comes from a stable model category structure on
differential graded $K(n)_{\ast}$-modules with weak equivalences 
the quasi-isomorphisms.  
The positive dimensional elements of $\pi_{\ast}^s$ act trivially on the 
homotopy categories in both cases. 
So the homotopy categories of the model categories of $K(n)$-modules 
and differential graded  $K(n)_{\ast}$-modules
are $\pi_{\ast}^s$-linearly equivalent.
However, the two model categories are not Quillen equivalent;
if they were Quillen equivalent, then the homotopy types of the function spaces
would agree \cite[Prop.\ 5.4]{DK}. But all function spaces of DG-modules
are products of Eilenberg-MacLane spaces, and this is not true for
$K(n)$-modules.
}
\end{remark}

\section{Margolis' uniqueness conjecture} \label{sec-margolis}

H.~R.~Margolis in `Spectra and the Steenrod algebra' introduced a set 
of axioms for a stable homotopy category \cite[Ch.\ 2 \S 1]{margolis}. 
The stable homotopy category of spectra satisfies the axioms, 
and Margolis conjectures \cite[Ch.\ 2, \S 1]{margolis} that this is 
the only model, i.e., that any category which satisfies the axioms 
is equivalent to the stable homotopy category. 
As part of the structure Margolis requires the subcategory of small objects
of a stable homotopy category to be equivalent to the Spanier-Whitehead
category of finite CW-complexes. 
So his uniqueness question really concerns possible `completions' 
of the category of finite spectra to a triangulated 
category with infinite coproducts.
Margolis shows \cite[Ch.\ 5 Thm.\ 19]{margolis} that
modulo phantom maps each model of his axioms is equivalent 
to the standard model.  
Moreover, in~\cite{christ-strick}, Christensen and Strickland show that
in any model the ideal of phantoms is equivalent to the phantoms in the 
standard model.

For objects $A$  and $X$ of  a triangulated category $\T$
we denote by $\T(A,X)_*$ the graded abelian homomorphism group
defined by $\T(A,X)_m=\T(A[m],X)$ for $m\in\mZ$, where $A[m]$ is the $m$-fold
shift of $A$. If $\T$ is a $\pi_*^s$-triangulated category, then the
groups $\T(A,X)_*$ form a graded $\pi_*^s$-module.
 
\begin{definition}\label{def-generator, small, margolis}
{\em An object $G$ of a triangulated category $\T$ 
is called a {\em weak generator} 
if it detects isomorphisms, i.e., a map $f:X\to Y$ is an isomorphism 
if and only if it induces an isomorphism between the graded abelian 
homomorphism groups $\T(G,X)_*$ and $\T(G,Y)_*$.  
An object $G$ of $\T$ is {\em small} if for any family of objects 
$\{A_i\}_{i\in I}$ whose coproduct exists the canonical map
\[ \bigoplus_{i\in I} \, \T(G, A_i) \ \varrow{1cm} \ 
\T(G, \coprod_{i\in I} A_i) \]
is an isomorphism.

A {\em stable homotopy category} in the sense of \cite[Ch.\ 2 \S 1]{margolis}
is a triangulated category $\mathscr{S}$ endowed with a symmetric monoidal, 
bi-exact smash product $\sm$ such that:
\begin{itemize}
\item $\mathscr{S}$ has infinite coproducts, 
\item the unit of the smash product is a small weak generator, and
\item there exists an exact and strong symmetric monoidal equivalence 
$R:\SW_{\mbox{\scriptsize f}}\to {\mathscr{S}}\rscript{small}$
between the Spanier-Whitehead category of finite CW-complexes 
(\cite{spanier-whitehead}, \cite[Ch.\ 1, \S 2]{margolis}) 
and the full subcategory of small objects in $\mathscr{S}$.
\end{itemize}}
\end{definition} 

The condition that $R$ is strong monoidal means that there
are coherently unital, associative, and commutative isomorphisms
between $R(A\sm B)$ and $R(A)\sm R(B)$ 
and between $R(S^0)$ and the unit of the smash product in $\mathscr S$.
Hence a stable homotopy category $\mathscr S$ becomes a $\pi_*^s$-triangulated
category as follows. The elements of $\pi_n^s$ are precisely
the maps from $S^n$ to $S^0$ in the Spanier-Whitehead category.
So given  $\alpha\in\pi_n^s=\SW(S^n,S^0)$ and $f:X\to Y$ in  $\mathscr S$
we can form $f\sm R(\alpha):X\sm R(S^n)\to Y\sm R(S^0)$.
Via the isomorphisms  $X\sm R(S^n)\iso X[n]\sm R(S^0)\iso X[n]$ and
$Y\sm R(S^0)\iso Y$ we obtain an element in $\mathscr S(X[n],Y)$
which we define to be $\alpha\cdot f$.
This $\pi_*^s$-action is unital, associative and bilinear because 
of the coherence conditions on the functor $R$.

As a consequence of our main theorem we can prove a special case of Margolis'
conjecture, namely we can show that a category satisfying his axioms is
equivalent to the homotopy category of spectra if it has {\em some}
underlying model category structure. 
Note that we do {\em not} ask for any kind of internal smash product 
on the model category which occurs in the following theorem. 

\begin{theorem} \label{thm-Margolis}
Suppose that $\mathscr{S}$ is a stable homotopy category
in the sense of \cite[Ch.\ 2 \S 1]{margolis} which supports a
$\pi_*^s$-linear equivalence with the homotopy category of some
stable model category.
Then $\mathscr{S}$ is equivalent to the stable homotopy category of spectra.
\end{theorem}
\begin{proof} Let $\C$ be a stable model category which
admits a $\pi_*^s$-linear equivalence $\Phi:\mathscr{S}\to\HoC$.
The image $X\in \HoC$ under $\Phi$ of the unit object of the smash product 
is a small weak generator for the homotopy category of $\C$.
Because the equivalence $\Phi$ is $\pi_*^s$-linear, 
$X$ satisfies condition (3) of our main theorem, and
so $\C$ is Quillen equivalent to the model category of spectra.
Thus the homotopy category of $\C$ and the category $\mathscr{S}$ are 
$\pi_*^s$-linearly equivalent to the ordinary stable homotopy category 
of spectra.
\end{proof}

\section{The $R$-local model structure for spectra}
\label{sec-local}

In this section we review the stable model category structure for spectra 
defined by Bousfield and Friedlander \cite[\S 2]{BF}
and  establish the $R$-local model structure
(Lemma \ref{lem-R-local model structure}).

A {\em spectrum} consists of a sequence $\{X_n\}_{n\geq 0}$ 
of pointed simplicial sets together
with maps $\sigma_n:S^1\sm X_n\to X_{n+1}$. A morphism $f:X\to Y$ of spectra
consists of maps of pointed simplicial sets $f_n:X_n\to Y_n$ for all
$n\geq 0$ such that $f_{n+1}\circ\sigma_n=\sigma_n\circ(1_{S^1}\sm f_n)$.
We denote the category of spectra by $\Spc$.  
A spectrum $X$ is an {\em $\Omega$-spectrum} if for all $n$ the simplicial set
$X_n$ is a Kan complex and the adjoint $X_n\to \Omega X_{n+1}$
of the structure map $\sigma_n$ is a weak homotopy equivalence.
The {\em sphere spectrum} $\mS$ is defined by 
$\mS_n=S^n=(S^1)^{\sm n}$, with structure maps the identity maps. 
The homotopy groups of a spectrum are defined by
\[ \pi_* X \ = \ \colim_i \pi_{i+*} |X_i| \]
A morphism of spectra is a {\em stable equivalence} if it
induces an isomorphism of homotopy groups.
A map $X\to Y$ of spectra is a {\em cofibration} if the map $X_0\to Y_0$
and the maps
\[ X_n \cup_{S^1\sm X_{n-1}} S^1\sm Y_{n-1} \ \varrow{1cm} \ Y_n \]
for $n\geq 1$ are cofibrations (i.e., injections) of simplicial sets.
A map of spectra is a {\em stable fibration} if it has the right
lifting property (see \cite[I p.\ 5.1]{Q}, \cite[3.12]{DS} 
or \cite[1.1.2]{hovey-book}) 
for the maps which are both cofibrations and stable equivalences.

Bousfield and Friedlander show in \cite[Thm.\ 2.3]{BF} that 
the stable equivalences, cofibrations and stable fibrations 
form a model category structure for spectra. 
A variation of their model category structure is the {\em $R$-local} 
model structure for $R$ a subring of the ring of rational numbers. 
The $R$-local model category structure is well known, 
but we were unable to find a reference in the literature.
A map of spectra is an {\em $R$-equivalence} if it induces an isomorphism
of homotopy groups after tensoring with $R$
and is an {\em $R$-fibration} if it has the 
right lifting property with respect to all maps that are cofibrations 
and $R$-equivalences.

\begin{lemma} \label{lem-R-local model structure}
Let $R$ be a subring of the ring of rational numbers.
Then the cofibrations, $R$-fibrations and $R$-equivalences make 
the category of  spectra into a model category,
referred to as the {\em $R$-local model category structure}. 
A spectrum is fibrant in the $R$-local model structure if and only if it is
an $\Omega$-spectrum with $R$-local homotopy groups.
\end{lemma}

We use `RLP' to abbreviate `right lifting property'. 
For one of the factorization axioms we need the small object argument
(see \cite[II 3.4 Remark]{Q} or \cite[7.12]{DS}) relative to a set 
$J=J\rscript{lv}\cup J\rscript{st} \cup J_R$ of maps of spectra which we now
define.
We denote by $\Delta[i]$, $\partial\Delta[i]$ and $\Lambda^{k}[i]$ 
respectively the simplicial
$i$-simplex, its boundary and its $k$-th horn 
(the union of all $(i-1)$-dimensional faces except the $k$-th one).
A subscript `$+$' denotes a disjoint basepoint. 
We denote by $F_n K$ the spectrum freely generated by a simplicial set $K$
in dimension $n$, i.e., $(F_n K)_j = S^{j-n}\sm K$ (where $S^m=*$ for $m<0$).
Hence $F_n K$ is a shift desuspension of the suspension spectrum of $K$.

First, $J\rscript{lv}$ is the set of maps of the form
\[ F_n \, \Lambda^{k}[i]_+ \ \varr{1cm}{} \ F_n \, \Delta[i]_+ \]
for $i,n\geq 0$ and $0 \leq k \leq i$.
Then $J\rscript{st}$ is the set of maps 
\[  Z_n^j\, \sm\, \partial\Delta[i]_+  \cup_{F_{n+j} S^j\, 
\sm\, \partial\Delta[i]_+} F_{n+j} S^j\, \sm\, \Delta[i]_+\ 
\varr{1cm}{ } \ Z_n^j\, \sm\, \Delta[i]_+\,  \]
for $i, j, n \geq 0$.  Here,
\[ Z_n^j \ = \ (F_{n+j} S^j\, \sm\, \Delta[1]_+) \cup_{F_{n+j} S^j \times 1} 
F_n S^0 \]
is the reduced mapping cylinder of the map 
$\lambda \mc F_{n+j} S^j \ \to \ F_n S^0$ 
which is the identity in spectrum levels above $n+j$. 
Note that $\lambda$ is a stable equivalence, but it is not a cofibration.
Since both source and target of $\lambda$ are cofibrant, the inclusion 
of the source into the mapping cylinder $Z^j_n$ is a cofibration. 
It is shown in \cite[Lemma A.3]{sch-theories} that the stable
fibrations of spectra are precisely the maps with the RLP
with respect to the set $J\rscript{lv}\cup J\rscript{st}$.

For every natural number $k$ we choose a finite pointed simplicial set $M_k$ 
which has the weak homotopy type of the mod-$k$ Moore space of dimension two. 
We let $J_R$ be the set of maps 
\[F_n \Sigma^m M_k \ \varrow{1cm} \  F_n \Sigma^m C(M_k) \] 
for all $m, n\geq 0$ and all natural numbers $k$ which are invertible in $R$, 
where $C(M_k)$ denotes the cone of the Moore space.

Now we prove a sequence of claims:
\renewcommand{\labelenumi}{(\alph{enumi})}
\begin{enumerate}
\item A map $X\to *$ has the RLP for the set 
$J=J\rscript{lv}\cup J\rscript{st} \cup J_R$ if and only if
$X$ is an $\Omega$-spectrum with $R$-local homotopy groups.
\item A map which is an $R$-equivalence and has the
RLP for $J$ is also an acyclic fibration in the stable model structure.
\item Every map can be factored as a composite $p\circ i$ where 
$p$ has the RLP for $J$ and $i$ is a cofibration and an $R$-equivalence and
is built from maps in $J$ by coproducts, pushouts and composition.
\item A map is an $R$-fibration if and only if it has the RLP for $J$.
\end{enumerate}

\smallskip

(a) 
The RLP for $(J\rscript{lv}\cup J\rscript{st})$ means that $X$ 
is stably fibrant, i.e., an $\Omega$-spectrum.
For $\Omega$-spectra the lifting property with respect to the map
$F_n \Sigma^m M_k \to F_n \Sigma^m C(M_k)$ 
means precisely that every element in the mod-$k$ homotopy group
\[ [F_n \Sigma^m M_k,X]^{\Ho(\Spc)} \ \iso \ \pi_0 \, \Omega^m \map(M_k, X_n) 
\ \iso \ \pi_{m+2-n} (X;\mZ/k) \] 
is trivial. 
Since this holds for all $m,n\geq 0$ and all $k$ which are invertible in $R$,
the map $X\to *$ has the RLP for $J$ if and only if  
$X$ is an $\Omega$-spectrum with $R$-local homotopy groups.

(b) Suppose $f:X\to Y$ is an $R$-equivalence and has the RLP for $J$.
Then $f$ is in particular a stable fibration and we denote its fiber by $F$.
There exists a long exact sequence connecting the homotopy groups of $F$,
$X$ and $Y$. Since $f$ is an $R$-equivalence, the localized homotopy 
groups $R\tensor\pi_*F$ of the fiber are trivial. 
As the base change of the map $f$, the map $F\to *$ also has the RLP for $J$. 
By (a), $F$ is an $\Omega$-spectrum whose homotopy groups are $R$-local.
Hence the homotopy groups of the fiber $F$ are trivial, so
the original map $f$ is also a stable equivalence.

(c) Every object occurring as the source of a map in $J$
is a suspension spectrum of a finite simplicial set, 
hence sequentially small in the sense of 
\cite[II 3.4 Remark]{Q} or \cite[Def.\ 7.14]{DS}.
Thus Quillen's {\em small object argument}
(see \cite[II 3.4 Remark]{Q} or \cite[7.12]{DS})
provides a factorization of a given map  as a composite $p\circ i$ 
where $i$ is built from maps in $J$ by coproducts, pushouts and composition, 
and where $p$ has the RLP for $J$.
Since every map in $J$ is a cofibration, so is $i$.
Cofibrations of spectra give rise to long exact sequences 
of homotopy groups, and homotopy groups of spectra commute 
with filtered colimits of cofibrations.
So to see that $i$ is an $R$-equivalence it suffices to check 
that the maps in $J$ are $R$-equivalences. 
The maps in $J\rscript{lv}$ are levelwise equivalences, 
the maps in $J\rscript{st}$ are stable equivalences, 
hence both are $R$-equivalences. 
Since the stable homotopy groups of the Moore space $M_k$ 
are $k$-power torsion, the maps in $J_R$ are also $R$-equivalences.

(d) We need to show that a map has the RLP for $J$ if and only if
it has the RLP for the (strictly bigger) class of maps $j$ 
which are cofibrations and $R$-equivalences.
This follows if any such $j$ is a retract of a map 
built from maps in $J$ by coproducts, pushouts and composition. 
We factor $j=p\circ i$ as in (c). Since $j$ and $i$ are $R$-equivalences,
so is $p$. Since $p$ also has the RLP for $J$, it is a stable acyclic
fibration by (b). So $p$ has the RLP for the cofibration $j$, hence $j$
is indeed a retract of $i$.

\smallskip

\begin{proof}[Proof of Lemma~\ref{lem-R-local model structure}] 
We verify the model category axioms as given in \cite[Def.\ 3.3]{DS}.
The category of spectra has all limits and colimits (MC1), 
the $R$-equivalences satisfy the 2-out-of-3 property (MC2) and the classes
of cofibrations, $R$-fibrations and $R$-equivalences 
are each closed under retracts (MC3).
By definition the $R$-fibrations have the RLP for maps which are both
cofibrations and $R$-equivalences.
Furthermore a map which is an $R$-equivalence and an $R$-fibration is
an acyclic fibration in the stable model structure by claim (b) above, 
so it has the RLP for cofibrations. This proves the lifting properties (MC4).
The stable model structure provides factorizations of maps as cofibrations 
followed by stable acyclic fibrations. 
Stable acyclic fibrations are in particular 
$R$-equivalences and $R$-fibrations, so this is also a factorization as a 
cofibration followed by an acyclic fibration in the $R$-local model structure.
The claims (c) and (d) provide the other factorization axiom (MC5).
\end{proof}

\begin{lemma} \label{R-local Lemma}
Let $\C$ be a stable model category, \mbox{$G:\C\to\Spc$} a functor 
with a left adjoint and $R$ a subring of the rational numbers. 
Then $G$ and its adjoint form a Quillen adjoint pair with respect to the
$R$-local model structure if and only if the following three conditions hold:
\renewcommand{\labelenumi}{(\roman{enumi})}
\begin{enumerate}
\item $G$ takes acyclic fibrations to level acyclic fibrations of spectra,
\item $G$ takes fibrant objects to $\Omega$-spectra with 
$R$-local homotopy groups and
\item $G$ takes fibrations between fibrant objects to level fibrations.
\end{enumerate}
\end{lemma}
\begin{proof} The `only if' part holds since the $R$-local acyclic fibrations
are level acyclic fibrations, the $R$-fibrant objects are the 
$\Omega$-spectra with $R$-local homotopy groups (claim (a) above),
and $R$-fibrations are in particular level fibrations.
For the converse suppose that $G$ satisfies conditions (i) to (iii).
We use a criterion of Dugger \cite[A.2]{dugger-simplicial}: 
in order to show that $G$ and its adjoint 
form a Quillen adjoint pair it suffices to show that $G$ 
preserves acyclic fibrations and it preserves fibrations between 
fibrant objects. The $R$-local acyclic fibrations are precisely 
the level acyclic fibrations, so  $G$ 
preserves acyclic fibrations by assumption (i).
We claim that every level fibration  $f:X\to Y$ between $\Omega$-spectra with 
$R$-local homotopy groups is an $R$-fibration. Given this,
$G$ preserves fibrations between fibrant objects by assumptions (ii) and (iii).

To prove the claim we choose a factorization $f=p\circ i$ with
$i:X\to Z$ a cofibration and $R$-equivalence and with $p:Z\to Y$
an $R$-fibration. Since $Y$ is $R$-fibrant, so is $Z$. Hence $i$
is an $R$-equivalence between $\Omega$-spectra with 
$R$-local homotopy groups, thus a level equivalence.
Hence $i$ is an acyclic cofibration in the strict model (or level)
model structure for spectra of \cite[2.2]{BF}, so that the
level fibration $f$ has the RLP for $i$. Hence $f$ is a retract of the
$R$-fibration $p$, and so it is itself an $R$-fibration. 
\end{proof}

\section{A universal property of the model category of spectra}
\label{spectra}

In this section we formulate a universal property which roughly says 
that the category of spectra is the `free stable model category on one object'.
The following theorem  associates to each cofibrant and fibrant object
$X$ of a stable model category $\C$ a Quillen adjoint functor pair such that
the left adjoint takes the sphere spectrum to $X$.
Moreover, this Quillen pair is essentially uniquely determined 
by the object $X$.
Theorem \ref{thm-recognition} gives conditions under which the adjoint pair
forms a Quillen equivalence.
We prove Theorem \ref{thm-universal} in the final section 
\ref{sec-constructions}. 

\begin{theorem} \label{thm-universal} {\bf (Universal property of spectra)} 
Let $\C$ be a stable model category and $X$ a cofibrant and fibrant object 
of $\C$.
\renewcommand{\labelenumi}{(\arabic{enumi})}
\begin{enumerate}
\item There exists a Quillen adjoint functor pair $X\sm -:\Spc\to\C$
and $\Hom(X,-):\C\to\Spc$ such that the left adjoint $X\sm -$ 
takes the sphere  spectrum, $\mS$, to $X$.
\item If $R$ is a subring of the rational numbers and the endomorphism group
$[X,X]^{\HoC}$ is an $R$-module, 
then any adjoint functor pair satisfying {\em (1)} is also a Quillen pair 
with respect to the $R$-local stable model structure for spectra.
\item If $\C$ is a {\em simplicial} model category, then the adjoint functors
 $X\sm -$ and $\Hom(X,-)$ of {\em (1)}
can be chosen as a {\em simplicial} Quillen adjoint functor pair.
\item Any two Quillen functor pairs satisfying {\em (1)} are related by a chain
of natural transformations which are weak equivalences on
cofibrant or fibrant objects respectively.
\end{enumerate}
\end{theorem}

Now we define the unit map
and deduce the $R$-local form of our main uniqueness theorem. 

\begin{definition} \label{def-unit}
{\em Let $X$ be a cofibrant and fibrant object of a stable model category $\C$.
Choose a Quillen adjoint pair  $X\sm -:\Spc\to\C$
and $\Hom(X,-):\C\to\Spc$ as in part (1) of Theorem \ref{thm-universal}.
The {\em unit map} of $X$ is the map of spectra
\[ \mS \ \varrow{1cm} \ \Hom(X,X) \]
which is adjoint to the isomorphism $X\sm\mS\iso X$.
By the uniqueness part (4) of Theorem \ref{thm-universal},
the spectrum $\Hom(X,X)$ is independent of the choice of Quillen pair 
up to stable equivalence of spectra under $\mS$.}
\end{definition}

\begin{theorem} \label{thm-recognition} Let $R$ be a subring of
the ring of rational numbers and let $\C$ be a stable model category. 
Then the following four conditions are equivalent:
\renewcommand{\labelenumi}{(\arabic{enumi})}
\begin{enumerate}
\item There is a chain of Quillen equivalences between $\C$ 
and the $R$-local stable model category of spectra.
\item There exists a $\pi_{\ast}^s$-linear equivalence between the homotopy 
category of $\C$ and the homotopy category of $R$-local spectra.
\item The homotopy category of $\C$ has a small weak generator $X$ for which 
$[X,X]_{\ast}^{\HoC}$ is freely generated as an 
$R\tensor \pi_{\ast}^s$-module by the identity map of $X$.
\item  The model category $\C$ has a cofibrant-fibrant 
small weak generator $X$ for which the  
groups  $[X,X]_{\ast}^{\HoC}$ are $R$-modules and the unit map  
$\mS\to \Hom(X,X)$ induces an isomorphism of homotopy groups after
tensoring with $R$.
\end{enumerate}
Furthermore, if $X$ is a cofibrant and fibrant object of $\C$ which satisfies
conditions {\em (3)} or {\em (4)}, then the functors $\Hom(X,-)$ and $X\sm -$ 
of Theorem \ref{thm-universal} {\em (1)} form a Quillen equivalence 
between $\C$ and the $R$-local model category of spectra.
\end{theorem}

\begin{remark}\label{rem-ss-class} 
{\em In \cite{ss-class} we associate to every object
of a stable model category an {\em endomorphism ring spectrum}. 
The spectrum $\Hom(X,X)$ given by Theorem \ref{thm-universal} (1) is stably 
equivalent to the underlying spectrum of the endomorphism ring spectrum. 
Moreover, the unit map as defined in \ref{def-unit} corresponds to
the unit map of ring spectra.
So condition (4) of the above theorem means that the endomorphism
ring spectrum of $X$ is stably equivalent, as a ring spectrum,
to the $R$-local sphere ring spectrum.
This expresses the equivalence of conditions (1) and (4) as a corollary 
of the more general classification result of \cite{ss-class}
for stable model categories with a small weak generator.
The special case in this paper, however, has a more direct proof.}
\end{remark}

\begin{proof}[Proof of Theorem \ref{thm-recognition}]
Every Quillen equivalence between stable model categories induces 
an equivalence of triangulated homotopy categories.
The derived functor of a left Quillen functor is also $\pi^s_*$-linear 
by Lemma \ref{left derived}; if it is an equivalence, then the inverse
equivalence is also $\pi^s_*$-linear. So condition (1) implies (2). 
Now assume (2) and let $X$ be a cofibrant and fibrant object
of $\HoC$ which in the homotopy category is isomorphic to the image of the 
localized sphere spectrum under any $\pi^s_*$-linear equivalence.
With this choice, condition (3) holds.

Given condition (3), we may assume that $X$ is cofibrant 
and fibrant and we choose a Quillen adjoint pair $X\sm -$ and $\Hom(X,-)$
as in part (1) of Theorem \ref{thm-universal}.
Since the group $[X,X]^{\HoC}$ is an $R$-module, 
the functors form a Quillen pair with respect to the $R$-local 
model structure for spectra by Theorem \ref{thm-universal} (2).
By Lemma \ref{left derived} the map
\[  X\sm^L - \ : \ [\mS,\mS]_*^{\Ho(\Spc_R)} \ \varrow{1cm} \ [X,X]_*^{\HoC} \]
induced by the left derived functor $X\sm^L-$ 
and the identification $X\sm^L\mS\iso X$ is $\pi_*^s$-linear
(note that the groups on the left hand side are taken in the $R$-local
homotopy category, so that $[\mS[n],\mS]^{\Ho(\Spc_R)}$ is isomorphic to
$R\tensor\pi_n^s$).
Source and target of this map are free $R\tensor\pi_*^s$-modules, 
and the generator Id$_{\mS}$ is taken to the generator Id$_X$.
Hence the map $X\sm^L -$ is an isomorphism.
For a fixed integer $n$, the derived adjunction and the identification
$X[n]\iso X\sm^L \mS[n]$ provide an isomorphism between $[X[n],X]^{\HoC}$ 
and $[\mS[n],\RHom(X,X)]^{\Ho(\Spc_R)}$ under which 
$X\sm^L-$ corresponds to $[\mS[n],\mS]^{\Ho(\Spc_R)} \to 
[\mS[n],\RHom(X,X)]^{\Ho(\Spc_R)}$ given by composition with the unit map.
For every spectrum $A$ the group $[\mS[n],A]^{\Ho(\Spc_R)}$ 
is naturally isomorphic to $R\tensor \pi_n A$, so this shows that the unit map
induces an isomorphism of homotopy groups
after tensoring with $R$, and condition (4) holds.

To conclude the proof we assume condition (4) and show that the 
Quillen functor pair $\Hom(X,-)$ and $X\sm -$ of
Theorem \ref{thm-universal} (1) is a Quillen equivalence.
Since the group $[X,X]^{\HoC}$ is an $R$-module, 
the functors form a Quillen pair with respect to the $R$-local 
model structure for spectra by Theorem \ref{thm-universal} (2).
So we show that the adjoint total derived functors 
$\RHom(X,-):\HoC\to\Ho(\Spc_R)$ and 
$X\sm^L -:\Ho(\Spc_R)\to \HoC$ are inverse equivalences of homotopy categories.
Note that the right derived functor $\RHom(X,-)$ is taken 
with respect to the $R$-local model structure on spectra.

For a fixed integer $n$, the derived adjunction and the identification 
$X\sm^L \mS[n]\iso X[n]$ provide a natural isomorphism 
\[ (*) \hspace*{1cm}  \pi_n \RHom(X,Y) \ \iso  \ 
[\mS[n],\RHom(X,Y)]^{\Ho(\Spc_R)} \ \iso \ [X[n],Y]^{\HoC} \ . \hspace*{1.5cm} \]
So the functor $\RHom(X,-)$ reflects isomorphisms because $X$ 
is a weak generator.
Hence it suffices to show that for every spectrum $A$ the unit 
of the adjunction of derived functors $A \to \RHom(X,X\sm^L A)$ 
is an isomorphism in the stable homotopy category.  
Consider the full subcategory $\T$ of the $R$-local stable homotopy category
with objects those spectra $A$ for which  $A \to \RHom(X,X\sm^L A)$ 
is an isomorphism.
Condition (4) says that the unit map $\mS\to \Hom(X,X)$ is an 
$R$-local equivalence, so $\T$ contains the (localized) sphere spectrum. 
Since the composite functor $\RHom(X,X\sm^L -)$ commutes with 
(de-)suspension and preserves distinguished triangles, 
$\T$ is a triangulated subcategory of the homotopy category of spectra.
As a left adjoint the functor $X\sm^L -$ preserves coproducts.
By formula $(*)$ above and since $X$ is small,
the natural map $\coprod_{I} \RHom(X,A_i) \to \RHom(X,\coprod_I A_i)$ is 
a $\pi_*$-isomorphism of spectra for any family of objects $A_i$ in $\HoC$. 
Hence the functor $\RHom(X,-)$ also preserves coproducts.
So $\T$ is a triangulated subcategory of the homotopy category 
of spectra which is also closed under coproducts 
and contains the localized sphere spectrum.
Thus, $\T$ is the whole $R$-local stable homotopy category, 
and this finishes the proof.
\end{proof}

\section{Construction of homomorphism spectra}
\label{sec-constructions}

In this last section we show that the derived functor of a
left Quillen functor is $\pi^s_*$-linear, and we prove 
Theorem \ref{thm-universal}. 

\begin{lemma}  \label{left derived}
Let $F:\C\to \D$ be the left adjoint of a Quillen adjoint
pair between stable model categories. Then the total
left derived functor $LF:\HoC\to\Ho(\D)$ is $\pi_*^s$-exact with respect
to the natural isomorphism $\tau:LF(X)\sm^L S^1\to LF(X\sm^L S^1)$
of \cite[5.6.2]{hovey-book}.
\end{lemma} 
\begin{proof} To simplify notation we abbreviate the derived functor
$LF$ to $L$ and drop the superscript $^{L}$ over the
smash product on the homotopy category level.
By \cite[5.7.3]{hovey-book}, the left derived functor $L$ is compatible
with the actions of the homotopy category of pointed simplicial sets --
Hovey summarizes this compatibility under the name of 
`$\Ho(\sset)$-module functor' \cite[4.1.7]{hovey-book}.
The isomorphism $\tau:L(X)\sm S^1\to L(X\sm S^1)$ is the special case 
$K=S^1$ of a natural isomorphism
\[ \tau_{X,K} \ : \ L(X) \ \sm K \ \varrow{1cm} \ L(X\sm K) \]
for a pointed simplicial set $K$ which is constructed in the proof 
of \cite[5.6.2]{hovey-book} 
(or rather its pointed analog in \cite[5.7.3]{hovey-book}).
It is important for us that the isomorphism $\tau$ is associative
(this is part of being a `$\Ho(\sset)$-module functor'), 
i.e., that the composite 
\[ L(A)\sm K\sm M \ \varr{1.8cm}{\tau_{A,K}\sm 1_M}  \  
L(A\sm K)\sm M   \ \varr{1.8cm}{\tau_{A\sm K,M}} \  L(A\sm K\sm M) \] 
is equal to $\tau_{A,K\sm M}$ (as before we suppress the implicit use of 
associativity isomorphisms such as $(A\sm K)\sm M\iso A\sm(K\sm M)$).
In particular the map 
$\tau_{X,S^n}:L(X)\sm S^n \to L(X\sm S^n)$ is equal to the $n$-fold
iterate of instances of $\tau_{-,S^1}$.

Now let $f:X\to Y$ be a morphism in the homotopy category of $\C$
and let $\alpha:S^{n+k}\to S^k$ represent a stable homotopy element.
We have to show that 
$\alpha\cdot L(f)=L(\alpha\cdot f)\circ\tau_{X,S^n}$ in the group
$[L(X)\sm S^n,L(Y)]^{\Ho(\D)}$.
By the definition of $\alpha\cdot L(f)$ this means proving
\begin{eqnarray} \label{eq1} L(f)\sm\alpha \ = \ 
(L(\alpha\cdot f)\circ\tau_{X,S^n})\sm 1_{S^k} 
\end{eqnarray}
in the group $[L(X)\sm S^{n+k},L(Y)\sm S^k]^{\Ho(\D)}$.
Since $\tau_{Y,S^{k}}:L(Y)\sm S^k\to L(Y\sm S^k)$ is an isomorphism
we may equivalently show equation (\ref{eq1}) after composition with
$\tau_{Y,S^{k}}$.
We note that 
\begin{eqnarray} \label{eq2}  \tau_{Y,S^k} \circ (L(f) \sm \alpha)
& = & L(f\sm \alpha) \circ \tau_{X,S^{n+k}} \\
\label{eq3} & = & L((\alpha\cdot f)\sm 1_{S^k}) \circ 
\tau_{X\sm S^n,S^{k}} \circ (\tau_{X,S^n} \sm 1_{S^k})\\
\label{eq4} & = &  \tau_{Y,S^k} \circ (L(\alpha\cdot f)\sm 1_{S^k}) \circ 
(\tau_{X,S^n} \sm 1_{S^k}) \\
 & = &  \tau_{Y,S^k} \circ ((L(\alpha\cdot f) \circ \tau_{X,S^n}) \sm 1_{S^k})
\ , \nonumber
\end{eqnarray}
which is what we had to show. 
Equations (\ref{eq2}) and (\ref{eq4}) use the naturality of $\tau$.
Equation (\ref{eq3}) uses the defining property of the morphism
$\alpha\cdot f$ and the associativity of $\tau$.
\end{proof}

Now we prove Theorem \ref{thm-universal}. We start with 

\begin{proof}[Proof of Theorem \ref{thm-universal} (2)]
By assumption the group $[X,X]^{\HoC}$ is a module over a subring $R$ 
of the ring of rational numbers.
Since $\Hom(X,-)$ is a right Quillen functor, it satisfies the
conditions of Lemma \ref{R-local Lemma} for $\mZ$. 
For fibrant $Y$, the $n$-th homotopy group of the $\Omega$-spectrum
$\Hom(X,Y)$ is isomorphic to the group 
$[\mS[n],\RHom(X,Y)]^{\Ho(\Spc)}$.
By the derived adjunction this group is isomorphic to the group 
$[X\sm^L\mS[n],Y]^{\HoC}\iso[X[n],Y]^{\HoC}$,
which is a module over the $R$-local endomorphism ring $[X,X]^{\HoC}$.
Hence the homotopy groups of the $\Omega$-spectrum $\Hom(X,Y)$ are $R$-local. 
Thus $\Hom(X,-)$ satisfies the conditions of Lemma \ref{R-local Lemma}
for $R$ and it is a right Quillen functor for the $R$-local model structure.
\end{proof}

Now we construct the adjoint functor pair $\Hom(X,-)$ and $X\sm -$ 
in the case of a {\em simplicial} stable model category.
This proves part (3) of Theorem \ref{thm-universal} and also
serves as a warm-up for the general construction which is very similar in 
spirit, but involves more technicalities.

\begin{construction} \label{con-simplicial}
{\em Let $\C$ be a {\em simplicial} stable model category and $X$ 
a cofibrant and fibrant object of $\C$.
We choose cofibrant and fibrant models $\omega^n X$ 
of the desuspensions of $X$ as follows. 
We set $\omega^0 X=X$ and inductively choose acyclic fibrations 
$\varphi_n:\omega^{n}X\to \Omega\,(\omega^{n-1}X)$ 
with $\omega^{n}X$ cofibrant.
We then define the functor $\Hom(X,-):\C \to \Spc$ by setting
\[ \Hom(X,Y)_n  \ = \ \map_{\C}(\omega^nX,Y)  \]
where `$\map_{\C}$' denotes  the simplicial mapping space. 
The spectrum structure maps are adjoint to the map
\[\begin{array}{rcl} \map_{\C}(\omega^{n-1}X,Y) & \ 
\varr{2cm}{\map_{\C}(\widetilde{\varphi_n},Y)} \ & 
\map_{\C}(\omega^{n}X\sm S^1, Y) \iso \Omega\, \map_{\C}(\omega^{n}X,Y) 
\end{array}\]
where $\widetilde{\varphi_n}$ is the adjoint of $\varphi_n$.

The functor $\Hom(X,-)$ has a left adjoint $X\sm -:\Spc\to \C$ 
defined as the coequalizer
\[ (\ast) \hspace{1.5cm} \bigvee_n \, \omega^n X \, \sm \, S^1  \sm \, A_{n-1} \
 \parallelarrows{1cm} \ \bigvee_n \, \omega^n X \, \sm \, A_n \ 
\varr{1cm} \ X \, \sm  \, A \ . \hspace{2cm} \]
The two maps in the coequalizer are induced by the structure maps 
of the spectrum $A$ and the maps 
$\widetilde{\varphi}_n: \omega^nX\sm S^1 \to \omega^{n-1} X$ respectively.
The various adjunctions provide bijections of morphism sets
\[ \C(X\sm\mS,W) \ \iso \ \Spc(\mS,\Hom(X,W)) \ \iso \ 
\sset(S^0,\Hom(X,W)_0) \ \iso \ {\C}(X,W) \]
natural in the $\C$-object $W$.
Hence the map $X\sm\mS\to X$ corresponding to the identity of $X$
in the case $W=X$ is an isomorphism; this shows that the left adjoint
takes the sphere spectrum to $X$.

Since $\omega^nX$ is cofibrant the functor $\map_{\C}(\omega^nX,-)$ 
takes fibrations (resp.\ acyclic fibrations) in $\C$ to fibrations 
(resp.\ acyclic fibrations)  of simplicial sets.
So the functor  $\Hom(X,-)$ takes  fibrations (resp.\ acyclic fibrations) 
in $\C$ to level fibrations (resp.\ level acyclic fibrations)
of spectra.
Since $\C$ is stable, $\widetilde{\varphi_n}$ is a weak equivalence
between cofibrant objects,
so for fibrant $Y$ the spectrum $\Hom(X,Y)$ is an $\Omega$-spectrum.
Hence $\Hom(X,-)$ satisfies the conditions of Lemma \ref{R-local Lemma} 
for $R=\mZ$, and so  $\Hom(X,-)$ and $X\sm -$ form a Quillen adjoint pair.
Since the functor $\Hom(X,-)$ is defined with the use of the
simplicial mapping space of $\C$, it comes with a natural, coherent
isomorphism $\Hom(X,Y^K)\iso \Hom(X,Y)^K$ for a simplicial set $K$.
So $\Hom(X,-)$ and its adjoint $X\sm -$ form a {\em simplicial}
Quillen functor pair which proves part (3) of Theorem \ref{thm-universal}.
}
\end{construction}

It remains to construct homomorphism spectra as in part (1) of Theorem 
\ref{thm-universal} for a general stable model category,
and prove the uniqueness part (4) of Theorem \ref{thm-universal}.
Readers who only work with simplicial model categories and have no need for
the uniqueness statement may safely ignore the rest of this paper.

To compensate for the lack of simplicial mapping spaces, 
we work with {\em cosimplicial frames}. 
The theory of `framings' of model categories goes back to Dwyer and Kan, 
who used the terminology {\em (co-)simplicial resolutions} \cite[4.3]{DK};
we mainly refer to Chapter 5 of Hovey's book \cite{hovey-book}
for the material about cosimplicial objects that we need.
If $K$ is a pointed simplicial set and $A$ a cosimplicial object of $\C$, 
then we denote by $A\sm K$ the coend
\cite[IX.6]{ML-working}
\[ A\sm K \  = \ \int^{n\in\Delta}  A^n \sm K_n\ , \]
which is an object of $\C$. Here $A^n\sm K_n$ denotes the coproduct of copies 
of $A^n$ indexed by the set $K_n$, modulo the copy of $A^n$ indexed by 
the basepoint of $K_n$.
Note that $A\sm\Delta[m]_+$ is naturally isomorphic to the object
of $m$-cosimplices of $A$; the object  $A\sm\partial\Delta[m]_+$
is also called the {\em $m$-th latching object} of $A$.
A cosimplicial map $A\to B$ is a {\em Reedy cofibration} if for all
$m\geq 0$ the map
\[ A\sm\Delta[m]_+ \cup_{A\sm\partial\Delta[m]_+} B\sm\partial\Delta[m]_+
\ \varrow{1cm} \ B\sm\Delta[m]_+ \]
is a cofibration in $\C$.
Cosimplicial objects in any pointed model category admit the    
{\em Reedy model structure} in which the weak equivalences are 
the cosimplicial maps which are levelwise weak equivalences 
and the cofibrations are the Reedy cofibrations. 
The Reedy fibrations are defined by the
right lifting property for Reedy acyclic cofibrations or equivalently
with the use of {\em matching objects}; see \cite[5.2.5]{hovey-book}
for details on the Reedy model structure.
If $A$ is a cosimplicial object and $Y$ is an object of $\C$, then there
is a simplicial set ${\C}(A,Y)$ of $\C$-morphisms defined by
${\C}(A,Y)_n={\C}(A^n,Y)$.
There is an adjunction bijection of pointed sets 
$\C(A\sm K,Y) \iso \sset(K,{\C}(A,Y))$.
If $A$ is a cosimplicial object, then the {\em suspension} of $A$
is the cosimplicial object $\Sigma A$ defined by
\[ (\Sigma A)^m \ = \ A\sm (S^1\sm\Delta[m]_+) \ . \]
Note that $\Sigma A$ and $A\sm S^1$ have different meanings:
$A\sm S^1$ is (naturally isomorphic to) the object of 0-cosimplices
of $\Sigma A$. 
There is a {\em loop} functor $\Omega$ for cosimplicial objects which is 
right adjoint to $\Sigma$; we do not use the precise form of $\Omega Y$ here.
For a cosimplicial object $A$ and an object $Y$ of $\C$ 
there is an adjunction isomorphism 
\[ {\C}(\Sigma A,Y) \ \iso \ \Omega \, {\C}(A,Y) \ . \]

A cosimplicial object in $\C$ {is \em homotopically constant} 
if each cosimplicial structure map is a weak equivalence in $\C$.
A {\em cosimplicial frame} (compare \cite[5.2.7]{hovey-book}) is a Reedy
cofibrant and homotopically constant cosimplicial object.
The following lemma collects from \cite[Ch.\ 5]{hovey-book}
those properties of cosimplicial frames which are relevant to our discussion.

\begin{lemma} \label{frame lemma} 
Let $\C$ be a pointed model category. 
\begin{enumerate}
\item The suspension functor for cosimplicial objects 
preserves Reedy cofibrations, Reedy acyclic cofibrations and
level equivalences between Reedy cofibrant objects. 
\item If $A$ is a cosimplicial frame, then so is $\Sigma A$.
\item If $A$ is a cosimplicial frame, then the functor 
${\C}(A,-)$ takes fibrations (resp.\ acyclic fibrations) in $\C$ to
fibrations (resp.\ acyclic fibrations) of simplicial sets.
\item If $Y$ is a fibrant object of $\C$, then the functor $\C(-,Y)$
takes level equivalences between Reedy cofibrant cosimplicial
objects to weak equivalences of simplicial sets.
\end{enumerate}
\end{lemma}
\begin{proof} 
(a) For a cosimplicial map $f\!:\!A\to B$ the map in $\C$
\[ (\Sigma A)\sm\Delta[m]_+ \cup_{(\Sigma A)\sm\partial\Delta[m]_+} 
(\Sigma B)\sm\partial\Delta[m]_+ \ \varrow{1cm} \ (\Sigma B)\sm\Delta[m]_+ \] 
is isomorphic to the pushout product $f\boxprod i$ \cite[4.2.1]{hovey-book}
of $f$ with the inclusion $i$ of $S^1\sm\partial\Delta[m]_+$ into 
$S^1\sm\Delta[m]_+$.
So if $f$ is a Reedy cofibration, then $f\boxprod i$ is a cofibration in
$\C$ by \cite[5.7.1]{hovey-book}; hence $\Sigma A\to\Sigma B$ is a Reedy
cofibration. 
In cosimplicial level $m$, the map $\Sigma f$ is given by the
map $f\sm (S^1\sm \Delta[m]_+)$. If $f$  is a Reedy acyclic 
cofibration, then $f\sm (S^1\sm \Delta[m]_+)$ is an acyclic cofibration
in $\C$ by \cite[5.7.1]{hovey-book}; hence $\Sigma f$ is also a level
equivalence. Suspension then preserves level equivalences between
Reedy cofibrant objects by Ken Brown's lemma \cite[1.1.12]{hovey-book}.

(b) If $A$ is a cosimplicial frame, then $\Sigma A$ is again Reedy cofibrant 
by part (a). A simplicial face map $d_i:\Delta[m-1]\to\Delta[m]$ induces 
an acyclic cofibration
\[ d_i^* \ : \ (\Sigma A)^{m-1} \ = \ A\sm(S^1\sm \Delta[m-1]_+)
\ \varrow{1cm} \ A\sm(S^1\sm \Delta[m]_+) \ = \  (\Sigma A)^{m} \]
by \cite[5.7.2]{hovey-book}, so $\Sigma A$ is also homotopically constant. 

(c) This is the pointed variant of \cite[5.4.4 (1)]{hovey-book}.

(d) If $A\to B$ is a Reedy acyclic cofibration, then for every 
cofibration of pointed simplicial sets $K\to L$ the map 
$A\sm L_{A\sm K} B\sm K\to B\sm L$ is an acyclic cofibration in $\C$ 
by \cite[5.7.1]{hovey-book}. By adjointness the induced map
${\C}(B,Y)\to{\C}(A,Y)$ is an acyclic fibration of simplicial sets.
By Ken Brown's Lemma \cite[1.1.12]{hovey-book}, the functor ${\C}(-,Y)$ 
thus takes level equivalences between Reedy cofibrant objects to 
weak equivalences of simplicial sets.
\end{proof}

The following lemma provides cosimplicial analogues
of the desuspensions $\omega^n X$ of Construction \ref{con-simplicial}.
 
\begin{lemma} \label{deloop}
Let $Y$ be a cosimplicial object in a stable model category $\C$ which is Reedy
fibrant and homotopically constant. Then there exists a
cosimplicial frame $X$ and a level equivalence $\Sigma X\to Y$
whose adjoint $X\to \Omega Y$ is a Reedy fibration which has the 
right lifting property for the map $*\to A$ for any cosimplicial frame $A$. 
\end{lemma}
\begin{proof}  Since $\C$ is stable there exists a cofibrant object $X^0$ 
of $\C$ such that the suspension of $X^0$ in the homotopy category of $\C$
is isomorphic to the object $Y^0$ of 0-cosimplices.
By \cite[4.5]{DK} or \cite[5.2.8]{hovey-book} there exists 
a cosimplicial frame $\bar X$ with ${\bar X}^0 = X^0$.
Since $\bar X$ is Reedy cofibrant, the map 
$d^0\amalg d^1:\bar{X}^0 \amalg \bar{X}^0 \to \bar{X}^1$ 
is a cofibration between cofibrant objects in $\C$;
since $\bar X$ is also homotopically constant, these maps express $\bar{X}^1$ 
as a cylinder object \cite[I 1.5 Def.\ 4]{Q} for $X^0$. 
The 0-cosimplices of $\Sigma \bar X$ are given by the quotient of the map  
$d^0\amalg d^1$, hence $(\Sigma \bar X)^0$ is a model 
for the suspension of $X^0$ in the homotopy category of $\C$.
Since $(\Sigma \bar X)^0$ is cofibrant and $Y^0$ is fibrant,
the isomorphism between them in the homotopy category can be realized 
by a weak equivalence $j^0:(\Sigma \bar X)^0\stackrel{\sim}{\to}Y^0$ 
in $\C$.
Since $Y$ is Reedy fibrant and homotopically constant, 
the map $Y\to cY^0$ is a Reedy acyclic fibration,
where $cY^0$ denotes the constant cosimplicial object. 
Since $\Sigma \bar X$ is Reedy cofibrant, the composite map
\[ \Sigma \bar X \ \varrow{1cm} \ c(\Sigma \bar X)^0 \ 
\varr{1cm}{cj^0} \ cY^0 \]
can be lifted to a map $j:\Sigma \bar X\to Y$.
The lift $j$ is a level equivalence since $j^0$ 
is an equivalence in $\C$
and both $\Sigma \bar X$ (by \ref{frame lemma} (b)) 
and $Y$ are homotopically constant.
The adjoint $\bar X\to \Omega Y$ of $j$ might not be a Reedy fibration,
but we can arrange for this by factoring it as a Reedy acyclic cofibration
$\bar X\to X$ followed by a Reedy fibration $\psi:X\to \Omega Y$, 
and replacing $j$ by the adjoint  
$\hat{\psi}:\Sigma X\to Y$ of the map $\psi$; 
by Lemma \ref{frame lemma} (a) the map $\Sigma \bar X\to\Sigma X$
is a level equivalence, hence so is $\hat\psi$.

Now suppose $A$ is a cosimplicial frame and $g:A\to \Omega Y$ is a cosimplicial
map with adjoint $\hat{g}:\Sigma A\to Y$. 
We want to construct a lifting, i.e., a map $A\to X$ whose
composite with $\psi:X\to\Omega Y$ is $g$.
We choose a cylinder object for $A$, i.e., a factorization
$A\Wedge A\to A\times I\to A$ of the fold map as a 
Reedy cofibration followed by a level equivalence. The suspension functor
preserves Reedy cofibrations and level equivalences between 
Reedy cofibrant objects by Lemma \ref{frame lemma} (a), 
so the suspended sequence 
$\Sigma A\Wedge \Sigma A\to \Sigma(A\times I)\to \Sigma A$ yields 
a cylinder object for $\Sigma A$.
In particular the 0-th level of  $\Sigma(A\times I)$ is a cylinder object 
for $(\Sigma A)^0=A\sm S^1$ in $\C$.

By \cite[6.1.1]{hovey-book} the suspension map 
$\Sigma:[A^0,X^0]\to [A^0\sm^L S^1,X^0\sm^L S^1]$
in the homotopy category of $\C$ can be constructed as follows. 
Given a $\C$-morphism $f^0:A^0\to X^0$,
one chooses an extension $f:A\to X$ to a cosimplicial map between
cosimplicial frames. The map $f\sm S^1:A\sm S^1\to X\sm S^1$ then 
represents the class $\Sigma [f^0]\in [A^0\sm^L S^1,X^0\sm^L S^1]$.
Composition with the 0-th level $\hat{\psi}^0:X\sm S^1\to Y^0$ 
of the level equivalence $\hat{\psi}:\Sigma X\to Y$ is a bijection from
$[A^0\sm^L S^1,X^0\sm^L S^1]$ to $[A^0\sm^L S^1,Y^0]$.
Since $\C$ is stable, the suspension map is bijective, which means that there
exists a cosimplicial map $f:A\to X$ such that the maps
$\hat{\psi}^0\circ (f\sm S^1)$ and $\hat{g}^0$ represent the same element in
$[A\sm S^1,Y^0]$.

The map $f$ need not be a lift of the original map $g$, but we can find
a lift in the homotopy class of $f$ as follows. 
Since $A\sm S^1$ is cofibrant and $Y^0$ is fibrant, there exists a 
homotopy $H_1:(\Sigma (A\times I))^0\to Y^0$ from 
$\hat{\psi}^0\circ (f\sm S^1)$ to $\hat{g}^0$. 
Evaluation at cosimplicial level zero is left adjoint
to the constant functor, so the homotopy $H_1$ is adjoint to a homotopy
$\hat{H}_1:\Sigma (A\times I)\to cY^0$ of cosimplicial objects.
Since $Y$ is Reedy fibrant and homotopically constant, the map
$Y\to cY^0$ is a Reedy acyclic fibration. So there exists a lifting
$H_2:\Sigma(A\times I)\to Y$ in the commutative square
\[\begin{diagram}
\node{\Sigma A\Wedge \Sigma A} \arrow{s,V} 
\arrow{e,t}{\hat{\psi}\circ(\Sigma f)\Wedge \hat{g}} \node{Y} \arrow{s,r,A}{\sim} \\
\node{\Sigma(A\times I)} \arrow{e,b}{\hat{H}_1} \node{cY^0}
\end{diagram} \]
which is a homotopy from $\hat{\psi}\circ (\Sigma f)$ to $\hat{g}$.
Taking adjoints gives a map $\hat{H}_2:A\times I\to \Omega Y$ which is
a homotopy from $\psi\circ f$ to $g$. 
Since $X\to \Omega Y$ is a Reedy fibration and the front
inclusion $i_0:A\to A\times I$ is a Reedy acyclic cofibration, 
we can choose a lifting $H_3:A\times I\to X$ in the commutative square
\[\begin{diagram}
\node{A} \arrow{s,l,V}{i_0} \arrow{e,t}{f} \node{X} \arrow{s,r,A}{\psi} \\
\node{A\times I} \arrow{e,b}{\hat{H}_2} \node{\Omega Y}
\end{diagram} \]
The end of the homotopy $H_3$, i.e., the composite map
$H_3\circ i_1:A\to X$, is then a lift of the original map $g:A\to \Omega Y$
since $\hat{H}_2 \circ i_1 = g$.
\end{proof}

\begin{construction} \label{con-general}
{\em Let $\C$ be a stable model category and $X$ 
a cofibrant and fibrant object of $\C$.
We define Reedy fibrant cosimplicial frames $\omega^n X$ as follows.
As in \cite[5.2.8]{hovey-book} we can choose a cosimplicial frame 
$\omega^0 X$ with $(\omega^0 X)^0=X$ and a Reedy acyclic fibration 
$\varphi_0:\omega^0 X\to cX$ which is the identity in dimension zero. 
Then $\omega^0 X$ is Reedy fibrant since $X$ is fibrant in $\C$. 
By Lemma \ref{deloop} we can inductively choose cosimplicial frames
$\omega^n X$ and level equivalences
$\hat{\varphi}_n:\Sigma (\omega^n X)\to \omega^{n-1} X$ whose
adjoints $\varphi_n:\omega^n X\to \Omega\,(\omega^{n-1}X)$ are Reedy fibrations
with the right lifting property for cosimplicial frames.
By Lemma \ref{frame lemma} (a), $\Sigma$ preserves Reedy acyclic
cofibrations, so $\Omega$ preserves Reedy fibrations. 
Hence $\Omega\,(\omega^{n-1} X)$ and thus $\omega^n X$ is Reedy fibrant.
We then define the functor $\Hom(X,-):\C \to \Spc$ by setting
\[ \Hom(X,Y)_n  \ = \ {\C}(\omega^nX,Y) \ . \]
The spectrum structure maps are adjoint to the map
\[\begin{array}{rcl} {\C}(\omega^{n-1}X,Y) & \ 
\varr{2cm}{{\C}(\hat{\varphi}_n,Y)} \ & 
{\C}(\Sigma\,(\omega^{n}X), Y) \iso \Omega\, {\C}(\omega^{n}X,Y) 
\end{array} \ . \]
The left adjoint $X\sm -:\Spc\to \C$ of $\Hom(X,-)$ is defined by the 
same coequalizer diagram $(*)$ as in 
Construction \ref{con-simplicial}, except that an expression like 
$\omega^n X \, \sm \, A_n$ now refers to the coend of a cosimplicial object 
with a simplicial set.
Also the isomorphism between $X\sm \mS$ and $X$ is obtained by the same
representability argument as in  \ref{con-simplicial}.

Since $\omega^nX$ is a cosimplicial frame, the functor ${\C}(\omega^nX,-)$ 
takes fibrations (resp.\ acyclic fibrations) in $\C$ 
to fibrations (resp.\ acyclic fibrations) 
of simplicial sets by Lemma \ref{frame lemma} (c).
So the functor  $\Hom(X,-)$ takes  fibrations (resp.\ acyclic fibrations) 
in $\C$ to level fibrations (resp.\ level acyclic fibrations) of spectra.
Since $\hat{\varphi}_n$ is a level equivalence between cosimplicial frames, 
Lemma \ref{frame lemma} (d) shows that the map ${\C}(\hat{\varphi}_n,Y)$ 
is a weak equivalence for fibrant $Y$; 
thus the spectrum $\Hom(X,Y)$ is an $\Omega$-spectrum for fibrant $Y$. 
So  $\Hom(X,-)$ and its adjoint form a Quillen pair by Lemma
\ref{R-local Lemma} for $R=\mZ$.
This proves part (1) of Theorem \ref{thm-universal}.}
\end{construction}

\begin{proof}[Proof of Theorem \ref{thm-universal} (4)]
Let $H:\Spc\to\C$ be any left Quillen functor with
an isomorphism $H(\mS)\iso X$, and let $G:\C\to\Spc$ be a right adjoint. 
We construct natural transformations $\Psi:\Hom(X,-)\to G$ and
$\Phi:H\to(X\sm -)$ where $\Hom(X,-)$ and $X\sm -$ are the Quillen pair
which were constructed in \ref{con-general}. 
Furthermore, $\Psi$ will be a stable equivalence of spectra for fibrant
objects of $\C$ and $\Phi$ will be a weak equivalence in $\C$ 
for every cofibrant spectrum.
So any two Quillen pairs as in Theorem \ref{thm-universal} (1) can be related
in this way via the pair $\Hom(X,-)$ and $X\sm -$ 
of Construction \ref{con-general}.

We denote by $F_n\Delta$ the cosimplicial spectrum given by
$(F_n\Delta)^m=F_n\Delta[m]_+$ and we denote by $H^{\bullet}$
the functor between cosimplicial objects obtained by applying
the left Quillen functor $H$ levelwise. The functor $H^{\bullet}$
is then a left Quillen functor with respect to the Reedy model structures
on cosimplicial spectra and cosimplicial objects of $\C$.
We inductively choose compatible maps 
$\psi_n:H^{\bullet}(F_n\Delta)\to \omega^n X$ of cosimplicial objects 
as follows.
Since $F_n\Delta$ is a cosimplicial frame,  
$H^{\bullet}(F_n\Delta)$ is a cosimplicial frame in $\C$.
The map $\varphi_0:\omega^0 X\to cX$ is a Reedy acyclic fibration, 
so the composite map 
\[ H^{\bullet}(F_0\Delta) \ \varrow{1cm} \ cH(F_0 S^0) \ 
\varr{1cm}{\iso} \ cX \]
admits a lift $\psi_0:H^{\bullet}(F_0\Delta)\to \omega^0 X$ which is
a level equivalence between cosimplicial frames.
The map $\varphi_n:\omega^nX\to \Omega\,(\omega^{n-1}X)$
has the right lifting property for cosimplicial frames, 
so we can inductively choose a lift 
$\psi_n:H^{\bullet}(F_n\Delta)\to \omega^n X$ of the composite map
\[ H^{\bullet}(F_n\Delta) \ \varrow{1cm} \ \Omega \, 
H^{\bullet}(F_{n-1}\Delta) 
\ \varr{1.5cm}{ \Omega(\psi_{n-1})} \ \Omega \, (\omega^{n-1}X) \ .\]
We show by induction that $\psi_n$ is a level equivalence.
The map $\psi_n\sm S^1$ is a weak equivalence in $\C$ since the other
three maps in the commutative square
\[\begin{diagram}
\node{\qquad H(F_nS^1) \iso H^{\bullet}(F_n\Delta)\sm S^1} \arrow{s} 
\arrow{e,t}{\psi_n\sm S^1}
\node{\omega^n X\sm S^1} \arrow{s,r}{\hat{\varphi}_n^0} \\
\node{\quad H(F_{n-1}S^0) \iso H^{\bullet}(F_{n-1}\Delta)^0}  
\arrow{e,b}{\psi_{n-1}^0}
\node{(\omega^{n-1} X)^0}
\end{diagram}\]
are. The map  $\psi_n\sm S^1$ is a model for the suspension of $(\psi_n)^0$. 
Since $\C$ is stable and  $(\psi_n)^0$ is a map between cofibrant objects,
 $(\psi_n)^0$ is a weak equivalence in $\C$. 
Since $H^{\bullet}(F_n\Delta)$ and $\omega^n X$ are homotopically constant,
the map $\psi_n:H^{\bullet}(F_n\Delta)\to\omega^n X$ is a level equivalence.

The adjunction provides a natural isomorphism of simplicial sets 
$G(Y)_n\iso {\C}(H^{\bullet}(F_n\Delta),Y)$ for every $n\geq 0$,
and we get a natural transformation 
\[ \Psi_n \ : \  \Hom(X,Y)_n \ = \ {\C}(\omega^nX,Y) \ 
\varr{2cm}{{\C}(\psi_n,Y)} \ {\C}(H^{\bullet}(F_n\Delta),Y) 
\ \iso \ G(Y)_n\ . \] 
By the way the maps $\psi_n$ were chosen,
the maps $\Psi_n$ together constitute a map of spectra 
$\Psi_Y:\Hom(X,Y)\to G(Y)$, natural in the $\C$-object $Y$. 
For fibrant objects $Y$, $\Psi_Y$ is a level equivalence, 
hence a stable equivalence, of spectra by Lemma \ref{frame lemma} (d)
since $\psi_n$ is a level equivalence between cosimplicial frames.

Now let $A$ be a spectrum. If we compose the adjoint 
$H(\Hom(X,X\sm A))\to X\sm A$ 
of the map $\Psi_{X\sm A}:\Hom(X,X\sm A)\to G(X\sm A)$ 
with $H(A) \to H(\Hom(X,X\sm A))$
coming from the adjunction unit, we obtain a natural transformation 
$\Phi_A:H(A)\to X \sm A$ between the left Quillen functors. 
The transformation $\Phi$ induces a natural transformation 
$L\Phi:LH\to X\sm^L -$ between the total left derived functors.
For any $Y$ in $\HoC$ the map  
$(L\Phi_A)^*:[X\sm^LA,Y]^{\HoC}\to [LH(A),Y]^{\HoC}$
is isomorphic to the bijection 
$(R\Psi_Y)_*:[A,\RHom(X,Y)]^{\Ho(\D)}\to[A,RG(Y)]^{\Ho(\D)}$.
Hence $L\Phi_A$ is an isomorphism in the homotopy category of $\C$
and so the map $\Phi_A$ is a weak equivalence in $\C$ 
for every cofibrant spectrum $A$.
\end{proof}

\end{document}